\documentclass[12pt,a4paper,reqno]{amsart}%{article}

\usepackage{amsmath}
\usepackage{amssymb}
\usepackage{amsfonts}

\advance\textheight by \topskip
\newtheorem{theorem}{Theorem}

\newtheorem{lemma}{Lemma}

\keywords{Binary search tree, 
 path length,   permutations, geometric distribution,
  $q$--analogues, harmonic numbers}
\subjclass{05A15, 05A30, 68P10}

\title[A $q$--analogue of the path length of binary 
search trees]
{A $q$--analogue of the path length of binary 
search trees}

\author{Helmut Prodinger}

\address{ Helmut Prodinger,
Centre for Applicable Analysis and Number Theory,
 Department of Mathematics,
University of the Witwatersrand, P.~O. Wits, 
2050 Johannesburg, South Africa, email:
{\tt helmut@gauss.cam.wits.ac.za}.
\newline
Homepage: {\tt http://www.wits.ac.za/helmut/index.htm}
}

\date{October 14, 1999}

\begin{document}

\begin{abstract}
A reformulation of the path length of binary search trees is
given in terms of permutations, allowing to extend the definition
to the instance of words, where the letters are obtained by
independent geometric random variables (with parameter $q$).
In this way, expressions for expectation and variance are
obtained which in the limit for $q\to1$ are the classical expressions. 
 
\end{abstract}

\maketitle

The path length $\rho(t)$ of a binary search tree
$t$ satisfies the recursion 
$\rho(t)=\rho(t_L)+\rho(t_R)+|t_L|+|t_R|$ where $t_L$ and
$t_R$   are the left resp. right subtree of the root. 
($|t| $ denotes the size of the tree $t$, i.~e. the number
of nodes.)

Binary search trees are obtained from permutations.
For some background see \cite{SeFl96, Knuth98, Mahmoud92}.
Our aim is to rewrite the definition of the path length in 
terms of permutations, since then we are able to obtain
$q$--analogues: This is done by considering words over
the alphabet $\{1,2,\dots\}$ instead, with probabilities
$p,pq,pq^2,\dots$, where $p+q=1$ (geometric probabilities).
In the limit $q\to1$, this model turns into the model of
random permutations, as equal letters appear with probability
0 and each relative ordering is equally likely.

For a permutation $\pi=\pi_1\dots\pi_n$   we define
$\rho(\pi)$ by 
\begin{equation*}   
\rho(\pi)=\big|\{(j,k)\mid 1\le j<k\le n,\
\pi_j=\min\{\pi_j,\dots,\pi_k\} \quad \text{or}\quad
\pi_k=\min\{\pi_j,\dots,\pi_k\}
\}\big|.
\end{equation*}

Then $\rho(\square)=0$ and, if $\pi=\sigma1\tau$, then
$\rho(\pi)=\rho(\sigma)+\rho(\tau)+|\sigma|+|\tau|$,
as pairs with the left coordinate in $\sigma$  and the
right coordinate in $\tau$ are definitely not counted. 

But this definition of $\pi$ can be taken as it is where
$\pi_1\dots\pi_n$ now denotes a {\it word\/} over the alphabet
$\{1,2,\dots\}$. This will be our starting point.

We want to point out that our previous paper \cite{prodinger99}
contains easier but related parameters.

In the sequel we want to compute the expectation and 
the variance of the parameter $\rho$, for random words
of length $n$. We define random variables

\begin{align*}
L_{jk}&=\begin{cases}1 & \text{if}\ \pi_j=\min\{\pi_j,\dots,\pi_k\}
 \\ 0 &\text{ otherwise},\end{cases}\\
R_{jk}&=\begin{cases}1 & \text{if}\ \pi_k=\min\{\pi_j,\dots,\pi_k\}
 \\ 0 &\text{ otherwise},\end{cases}\\
B_{jk}&=L_{jk}\cdot R_{jk},\\
N_{jk}&=(1-L_{jk})\cdot (1-R_{jk}).\\
\end{align*}

(The letters $L,R,B,N$ are chosen to indicate
{\it left, right, both, not}.)

Then the parameter $\rho$ may be described as
\begin{equation*}   
\rho=\sum_{1\le j< k\le n}\Big[L_{jk}+R_{jk}-B_{jk}\Big].
\end{equation*}

Now we can introduce  the generating function 
\begin{align*}
f(v)=\Big(\frac pq\Big)^{n}
\sum_{i_1,\dots,i_n\ge1}
q^{i_1+\dots+i_n}
\prod_{1\le j<k\le n}
\Big[L_{jk}v+R_{jk}v-B_{jk}v+N_{jk}
\Big];
\end{align*}
the coefficient of $v^k$ in $f(v)$ is the probability that parameter
$\rho$ has value $k$, assuming random words of length
$n$.

As always,  the expected value is  obtained via
$\mathbb{E}=f{'}(1)$;

\begin{align*}
\mathbb{E}&=\Big(\frac pq\Big)^{n}
\sum_{i_1,\dots,i_n\ge1}
q^{i_1+\dots+i_n}
\sum_{1\le j<k\le n}
\big(L_{jk}+R_{jk}-B_{jk}\big)
\\
&=\sum_{1\le j<k\le n}
\Big(\frac pq\Big)^{k+1-j}\Big[
\sum_{L_{jk}=1}q^{i_j+\dots+i_k}
+\sum_{R_{jk}=1}q^{i_j+\dots+i_k}
-\sum_{ B_{jk}=1}q^{i_j+\dots+i_k}\Big]
\\
&=\sum_{1\le j<k\le n}
\Big(\frac pq\Big)^{k+1-j}\Big[
2\sum_{i\ge1}q^{i(k+1-j)}\frac{1}{p^{k-j}}
-\sum_{i\ge1}q^{i(k+1-j)}\frac{1}{p^{k-j-1}}
\Big]
\\
&=(2p-p^2)\sum_{1\le j<k\le n}
\Big(\frac 1q\Big)^{k+1-j}
\sum_{i\ge1}q^{i(k+1-j)}\\
&=p(2-p)\sum_{1\le j<k\le n}
\frac{1}{1-q^{k+1-j}}\\
&=p(2-p)\sum_{2\le i\le n}\frac{n+1-i}{1-q^{i}}\\
&=p(2-p)\sum_{1\le i\le n}\frac{n+1-i}{1-q^{i}}-n(2-p).
\end{align*}

The terms that would survive the limit $q\to1$ are

\begin{equation*}   
2p\sum_{1\le i\le n}(n+1-i)
\frac{1}{1-q^{i}}-2n,
\end{equation*}
and the limit is
\begin{equation*}   
\lim_{q\to1}\mathbb{E}
=2\sum_{1\le i\le n}\frac{n+1-i}{i}-2n
=2(n+1)H_n-4n,
\end{equation*}
as is of course well known. 

Now we turn to the variance, and this is much harder, 
since we must first compute 
the second factorial moment, which is  obtained by
a second derivative;

\begin{align*}
\mathbb{E}^{\underline{2}}&=\Big(\frac pq\Big)^{n}
\sum_{i_1,\dots,i_n\ge1}
q^{i_1+\dots+i_n}\times\\ &\times
\sum_{1\le j<k\le n, 1\le l<m\le n, (j,k)\neq(l,m)}
\big(L_{jk}+R_{jk}-B_{jk}\big)
\big(L_{lm}+R_{lm}-B_{lm}\big)\\
&=\Big(\frac pq\Big)^{n}
\sum_{i_1,\dots,i_n\ge1}
q^{i_1+\dots+i_n}\times\\ &\times
\sum_{1\le j<k\le n, 1\le l<m\le n, (j,k)\neq(l,m)}
\big(2L_{jk}L_{lm}+2L_{jk}R_{lm}-4L_{jk}B_{lm}
+B_{jk}B_{lm}\big)\\
&=2\Xi^{\text{LL}}+2\Xi^{\text{LR}}-4\Xi^{\text{LB}}+\Xi^{\text{BB}}
\end{align*}
(using several symmetries).

The range $\Lambda=\{1\le j<k\le n, 1\le l<m\le n, (j,k)\neq(l,m)\}$
must be split into the following 12 disjoint subranges:

\begin{align*}
\Lambda_{1}&=\{1\le j<k<l<m\le n\},\\
\Lambda_{2}&=\{1\le j<l<m<k\le n\},\\
\Lambda_{3}&=\{1\le j<l<k<m\le n\},\\
\Lambda_{4}&=\{1\le j<k=l<m\le n\},\\
\Lambda_{5}&=\{1\le j<l<m=k\le n\},\\
\Lambda_{6}&=\{1\le j=l<k<m\le n\},\\
\Lambda_{7}&=\{1\le l<m<j<k\le n\},\\
\Lambda_{8}&=\{1\le l<j<k<m\le n\},\\
\Lambda_{9}&=\{1\le l<j<m<k\le n\},\\
\Lambda_{10}&=\{1\le l<m=j<k\le n\},\\
\Lambda_{11}&=\{1\le l<j<k=m\le n\},\\
\Lambda_{12}&=\{1\le l=j<m<k\le n\}.\\
\end{align*}

And we will have contributions $\Theta_i^{\text{LL}}$,
$\Theta_i^{\text{LR}}$,
$\Theta_i^{\text{LB}}$,
$\Theta_i^{\text{BB}}$, 
to $\Xi^{\text{LL}}$,
$\Xi^{\text{LR}}$,
$\Xi^{\text{LB}}$,
$\Xi^{\text{BB}}$,
for $i=1,\dots,12$, according to the 12
ranges $\Lambda_i$.

Therefore we must compute 48 (not necessarily) different 
contributions.

For convenience, we state them as a lemma.

\begin{lemma} The contributions
$\Theta_i^{\text{LL}}$,
$\Theta_i^{\text{LR}}$,
$\Theta_i^{\text{LB}}$,
$\Theta_i^{\text{BB}}$, for $i=1,\dots,12$, are given by

\begin{align*}
\Theta_{1}^{\text{LL}}&=\Theta_{7}^{\text{LL}}=
p^2\sum_{1\le j<k<l<m\le n}\frac{1}{1-q^{k+1-j}}
\frac{1}{1-q^{m+1-l}},\\
\Theta_{2}^{\text{LL}}&=\Theta_{8}^{\text{LL}}=
p^2\sum_{1\le j<l<m<k\le n}\frac{1}{1-q^{k+1-j}}
\frac{1}{1-q^{m+1-l}},\\
\Theta_{3}^{\text{LL}}&=\Theta_{9}^{\text{LL}}=
p^2\sum_{1\le j<l<k<m\le n}\frac{1}{1-q^{m+1-j}}
\frac{1}{1-q^{m+1-l}},\\
\Theta_{4}^{\text{LL}}&=\Theta_{10}^{\text{LL}}=
p^2\sum_{1\le j<k=l<m\le n}\frac{1}{1-q^{m+1-j}}
\frac{1}{1-q^{m+1-k}},\\
\Theta_{5}^{\text{LL}}&=\Theta_{11}^{\text{LL}}=
p^2\sum_{1\le j<l<m=k\le n}
\frac{1}{1-q^{m+1-j}}
\frac{1}{1-q^{m+1-l}},\\
\Theta_{6}^{\text{LL}}&=\Theta_{12}^{\text{LL}}=
p\sum_{1\le j=l<k<m\le n}
\frac{1}{1-q^{m+1-j}};
\end{align*}

\begin{align*}
\Theta_{1}^{\text{LR}}&=\Theta_{7}^{\text{LR}}=
p^2\sum_{1\le j<k<l<m\le n}\frac{1}{1-q^{k+1-j}}
\frac{1}{1-q^{m+1-l}},\\
\Theta_{2}^{\text{LR}}&=\Theta_{8}^{\text{LR}}=
p^2\sum_{1\le j<l<m<k\le n}\frac{1}{1-q^{k+1-j}}
\frac{1}{1-q^{m+1-l}},\\
\Theta_{3}^{\text{LR}}&=\Theta_{9}^{\text{LR}}=
p^2\sum_{1\le j<l<k<m\le n}\bigg[
\frac{1}{1-q^{m+1-j}}\frac{1}{1-q^{m+1-l}}+
\frac{1}{1-q^{m+1-j}}\frac{1}{1-q^{k+1-j}}-
\frac{1}{1-q^{m+1-j}}
\bigg],\\
\Theta_{4}^{\text{LR}}&
=p^2\sum_{1\le j<k=l<m\le n}\bigg[
\frac{1}{1-q^{m+1-j}}\frac{1}{1-q^{m+1-k}}+
\frac{1}{1-q^{m+1-j}}\frac{1}{1-q^{k+1-j}}-
\frac{1}{1-q^{m+1-j}}
\bigg],\\
\Theta_{5}^{\text{LR}}&=
p^2\sum_{1\le j<l<m\le n}
\frac{1}{1-q^{m+1-j}}
\frac{1}{1-q^{m+1-l}},\\
\Theta_{6}^{\text{LR}}&
=p^2\sum_{1\le j<k<m\le n}
\frac{1}{1-q^{m+1-j}}\frac{1}{1-q^{k+1-j}},\\
\Theta_{10}^{\text{LR}}&
=p\sum_{1\le l<m=j<k\le n}
\frac{1}{1-q^{k+1-l}},\\
\Theta_{11}^{\text{LR}}&=\Theta_{12}^{\text{LR}}
=p^2\sum_{1\le l=j<m<k\le n}
\frac{1}{1-q^{k+1-l}};\\
\end{align*}

\begin{align*}
\Theta_{1}^{\text{LB}}&=\Theta_{7}^{\text{LB}}=
p^3\sum_{1\le j<k<l<m\le n}\frac{1}{1-q^{k+1-j}}
\frac{1}{1-q^{m+1-l}},\\
\Theta_{2}^{\text{LB}}&=\Theta_{8}^{\text{LB}}=
p^3\sum_{1\le j<l<m<k\le n}\frac{1}{1-q^{k+1-j}}
\frac{1}{1-q^{m+1-l}},\\
\Theta_{3}^{\text{LB}}&=\Theta_{9}^{\text{LB}}=
p^3\sum_{1\le j<l<k<m\le n}
\frac{1}{1-q^{m+1-j}}\frac{1}{1-q^{m+1-l}}
,\\
\Theta_{4}^{\text{LB}}&
=p^3\sum_{1\le j<k=l<m\le n}
\frac{1}{1-q^{m+1-j}}\frac{1}{1-q^{m+1-k}}
,\\
\Theta_{5}^{\text{LB}}&=
p^3\sum_{1\le j<l<m\le n}
\frac{1}{1-q^{m+1-j}}
\frac{1}{1-q^{m+1-l}},\\
\Theta_{6}^{\text{LB}}&
=p^2\sum_{1\le j<k<m\le n}
\frac{1}{1-q^{m+1-j}},\\
\Theta_{10}^{\text{LB}}&
=p^2\sum_{1\le l<m=j<k\le n}
\frac{1}{1-q^{k+1-l}},\\
\Theta_{11}^{\text{LB}}&=
p^3\sum_{1\le l<j<k=m\le n}
\frac{1}{1-q^{k+1-l}},\\
\Theta_{12}^{\text{LB}}
&=p^2\sum_{1\le l=j<m<k\le n}
\frac{1}{1-q^{k+1-l}};\\
\end{align*}

\begin{align*}
\Theta_{1}^{\text{BB}}&=\Theta_{7}^{\text{BB}}=
p^4\sum_{1\le j<k<l<m\le n}\frac{1}{1-q^{k+1-j}}
\frac{1}{1-q^{m+1-l}},\\
\Theta_{2}^{\text{BB}}&=\Theta_{8}^{\text{BB}}=
p^4\sum_{1\le j<l<m<k\le n}\frac{1}{1-q^{k+1-j}}
\frac{1}{1-q^{m+1-l}},\\
\Theta_{3}^{\text{BB}}&=\Theta_{9}^{\text{BB}}=
p^4\sum_{1\le j<l<k<m\le n}\frac{1}{1-q^{m+1-j}}
,\\
\Theta_{4}^{\text{BB}}&=\Theta_{10}^{\text{BB}}=
p^3\sum_{1\le j<k=l<m\le n}\frac{1}{1-q^{m+1-j}}
,\\
\Theta_{5}^{\text{BB}}&=\Theta_{11}^{\text{BB}}=
p^3\sum_{1\le j<l<m\le n}
\frac{1}{1-q^{m+1-j}}
\frac{1}{1-q^{m+1-l}},\\
\Theta_{6}^{\text{BB}}&=\Theta_{12}^{\text{BB}}=
p^3\sum_{1\le j=l<k<m\le n}
\frac{1}{1-q^{m+1-j}}.
\end{align*}
\end{lemma}
\begin{proof}
The computations are as (or slightly more complicated than)
the one for the expected value. We don't give more details. 
\end{proof}

We  simplify those sums and
write $a_i=\frac1{1-q^i}$ for convenience.

\begin{lemma}
\begin{align*}
\Xi^{\text{LL}}&=2p^2\sum_{2\le i,j\le n-2, i+j\le n}a_ia_j\binom{n+2-i-j}{2}\\
&+2p^2\sum_{2\le i<j\le n}a_ia_j(n+1-j)(j-i-1)\\
&+2p^2\sum_{3\le i<j\le n}a_ia_j(n+1-j)(i-2)\\
&+4p^2\sum_{2\le i<j\le n}a_ia_j(n+1-j)\\
&+2p\sum_{3\le i\le n}a_i(i-2)(n+1-i),
\end{align*}

\begin{align*}
\Xi^{\text{LR}}&=2p^2\sum_{2\le i,j\le n-2, i+j\le n}a_ia_j\binom{n+2-i-j}{2}\\
&+2p^2\sum_{2\le i<j\le n}a_ia_j(n+1-j)(j-i-1)\\
&+4p^2\sum_{3\le i<j\le n}a_ia_j(n+1-j)(i-2)\\
&-2p^2\sum_{4\le i\le n}a_i\binom {i-2}2(n+1-i)\\
&+4p^2\sum_{2\le i<j\le n}a_ia_j(n+1-j)\\
&+p^2\sum_{3\le i\le n}a_i({i-2})(n+1-i)\\
&+p\sum_{3\le i\le n}a_i(i-2)(n+1-i),
\end{align*}

\begin{align*}
\Xi^{\text{LB}}&=2p^3\sum_{2\le i,j\le n-2, i+j\le n}a_ia_j\binom{n+2-i-j}{2}\\
&+2p^3\sum_{2\le i<j\le n}a_ia_j(n+1-j)(j-i-1)\\
&+2p^3\sum_{3\le i<j\le n}a_ia_j(n+1-j)(i-2)\\
&+2p^3\sum_{2\le i<j\le n}a_ia_j(n+1-j)\\
&+(3p^2+p^3)\sum_{3\le i\le n}a_i(i-2)(n+1-i),
\end{align*}

\begin{align*}
\Xi^{\text{BB}}&=2p^4\sum_{2\le i,j\le n-2, i+j\le n}a_ia_j\binom{n+2-i-j}{2}\\
&+2p^4\sum_{2\le i<j\le n}a_ia_j(n+1-j)(j-i-1)\\
&+2p^3\sum_{2\le i<j\le n}a_ia_j(n+1-j)\\
&+2p^4\sum_{4\le i\le n}a_i\binom{i-2}{2}(n+1-i)\\
&+4p^3\sum_{3\le i\le n}a_i({i-2})(n+1-i).
\end{align*}

\qed

\end{lemma}

The variance is given by 
\begin{equation*}   
\mathbb{V}=2\Xi^{\text{LL}}+2\Xi^{\text{LR}}-4\Xi^{\text{LB}}+\Xi^{\text{BB}}
+\mathbb{E}-\big(\mathbb{E}\big)^2.
\end{equation*}
In order to simplify this expression, we note the following
formul{\ae}:

\begin{lemma}
\begin{align*}
p^2\sum_{2\le i,j\le n-2, i+j\le n}&a_ia_j\binom{n+2-i-j}{2}\\
&=2p^2\sum_{1\le i<j\le  n}a_ia_j\binom{n+2-j}{2}
-p^2\sum_{1\le i\le  n}a_i\binom{n+2-i}{2}(i-1)\\&
-2p\sum_{1\le i \le n}a_i\binom{n+1-i}2+\binom n2.
\end{align*}
\end{lemma}
\begin{proof}
Note that
\begin{equation*}
\frac{1}{1-q^i}\frac{1}{1-q^j}=
\frac{1}{1-q^{i+j}}\bigg(
\frac{1}{1-q^i}+
\frac{1}{1-q^j}-1
\bigg)
\end{equation*}
and do some trivial rearrangements. 
\end{proof}

\begin{lemma}
\begin{align*}
\Big(\sum_{1\le i\le n}a_i(n+1-i)\Big)^2=
2\sum_{1\le i<j\le n}a_ia_j(n+1-i)(n+1-j)
+\sum_{1\le i\le n}a_i^2(n+1-i)^2.
\end{align*}

\end{lemma}

\begin{proof}Obvious.\end{proof}

Using these lemmata and numerous simplifications that
were partially supported by Maple, we can state our main 
result:

\begin{theorem}
The expectation and the variance of the $q$--ified
path length in words of length $n$, generated by $n$ 
independent geometric random variables are given by

\begin{align*}
\mathbb{E}=p(2-p)\sum_{1\le i\le n}
\frac{n+1-i}{1-q^{i}}-n(2-p)
\end{align*}
and

\begin{align*}
\mathbb{V}&=2p^2\sum _{1\le i<j\le n}\frac{(n+1-j)(4i+p(5-4i) )}
{(1-q^i)(1-q^j)}\\
&-{p}^{2}(2-p )^{2}\sum _{1\le i\le n}
\frac{(n+1-i)^{2}}{(1-q^i)^2}\\
&+p\sum_{1\le i\le n}\frac{n+1-i}{1-q^i}
\Big(6i-2 +p(-4ni+4n-19i+3{i}^{2}+7 )\\&\qquad\qquad+4p^2 (ni-n+3i
-1-{i}^{2} )+ {p}^{3}(-ni+n+2{i}^{2}-8i+8
 )\Big)\\
&+5pn-3{p}^{2}n-2{p}^{3}n.
\end{align*}

\qed 

\end{theorem}

The terms in the variance that would survive the limit
$q\to1$  are these:
\begin{align*}
8p^2\sum _{1\le i<j\le n}\frac{(n+1-j)i}{{(1-q^i)(1-q^j)}}
-4{p}^{2}\sum _{1\le i\le n}\frac{
(n+1-i)^{2}}{{(1-q^i)^2}}+2p\sum_{1\le i\le n}\frac{(n+1-i)
(3i-1)}{1-q^i}.
\end{align*}

The limit is

\begin{align*}
\lim_{q\to1}\mathbb{V}&=8\sum _{1\le j\le n}\frac{(n+1-j)(j-1)}{j}-4\sum _{1\le i\le n}
\frac{(n+1-i)^{2}}{i^2}+2\sum_{1\le i\le n}\frac{n+1-i}{i}
(3i-1)\\
&=4(n+1)n-8(n+1)H_n+8n-4(n+1)^2H_n^{(2)}+8(n+1)H_n-4n\\
&\qquad\qquad\qquad\qquad\qquad\qquad\qquad+3n(n+1)-2(n+1)H_n+2n\\
&=7n^2-4(n+1)^2H_n^{(2)}-2(n+1)H_n+13n.
\end{align*}

This is (of course!) the variance in the classical case.

\bibliographystyle{plain}

%\bibliography{pro_bib}

\end{document}